\documentclass[12pt]{amsart}    
\usepackage{amscd}

\def\cal#1{\mathcal{#1}}

\def\NZQ{\Bbb}               
\def\NN{{\NZQ N}}

\def\ZZ{{\NZQ Z}}

\def\frk{\frak}               

\def\mm{{\frk m}}

\def\opn#1#2{\def#1{\operatorname{#2}}} 
\opn\chara{char}
\opn\length{\ell}
\opn\pd{pd}
\opn\rk{rk}
\opn\projdim{proj\,dim}
\opn\rank{rank}
\opn\depth{depth}
\opn\grade{grade}
\opn\height{height}
\opn\embdim{emb\,dim}
\opn\codim{codim}

\opn\Tr{Tr}
\opn\bigrank{big\,rank}
\opn\superheight{superheight}\opn\lcm{lcm}
\opn\trdeg{tr\,deg}%
\opn\reg{reg}
\opn\lreg{lreg}
\opn\div{div}
\opn\Div{Div}
\opn\cl{cl}
\opn\Cl{Cl}
\opn\Spec{Spec}
\opn\Supp{Supp}
\opn\supp{supp}
\opn\Sing{Sing}
\opn\Ass{Ass}
\opn\Ann{Ann}
\opn\Rad{Rad}
\opn\Soc{Soc}
\opn\Ker{Ker}
\opn\Coker{Coker}
\opn\Im{Im}
\opn\Hom{Hom}
\opn\Tor{Tor}
\opn\Ext{Ext}
\opn\End{End}
\opn\Aut{Aut}
\opn\id{id}

\opn\nat{nat}
\opn\pff{pf}
\opn\Pf{Pf}
\opn\GL{GL}
\opn\SL{SL}
\opn\mod{mod}
\opn\ord{ord}
\opn\aff{aff}
\opn\con{conv}
\opn\relint{relint}
\opn\st{st}
\opn\lk{lk}
\opn\cn{cn}
\opn\core{core}
\opn\vol{vol}
\opn\link{link}
\opn\star{star}
\opn\gr{gr}

\def\poly#1#2#3{#1[#2_1,\dots,#2_{#3}]}
\def\pot#1#2{#1[\kern-0.28ex[#2]\kern-0.28ex]}

\opn\dirlim{\underrightarrow{\lim}}
\opn\inivlim{\underleftarrow{\lim}}
\let\tensor=\otimes
\let\iso=\cong

\let\Dirsum=\bigoplus

\let\mcone= * 
\let\To=\longrightarrow
\def\Implies{\ifmmode\Longrightarrow \else
     \unskip${}\Longrightarrow{}$\ignorespaces\fi}
\def\implies{\ifmmode\Rightarrow \else
     \unskip${}\Rightarrow{}$\ignorespaces\fi}
\def\iff{\ifmmode\Longleftrightarrow \else
     \unskip${}\Longleftrightarrow{}$\ignorespaces\fi}

\let\:=\colon
\newtheorem{Theorem}{Theorem} 
\newtheorem{Lemma}{Lemma} 
\newtheorem{Corollary}{Corollary} 
\newtheorem{Proposition}{Proposition} 
\newtheorem{Remark}{Remark}

\newtheorem{Example}{Example} 

\let\epsilon\varepsilon
\let\kappa=\varkappa
\opn\inii{in}
\opn\inim{inm}
\opn\set{set}
\def\pnt{{\raise0.5mm\hbox{\large\bf.}}}

\begin{document}

\title{
Combinatorial characterizations of 
generalized Cohen-Macaulay monomial ideals
}
\author{Yukihide Takayama}
\address{Yukihide Takayama, Department of Mathematical
Sciences, Ritsumeikan University, 
1-1-1 Nojihigashi, Kusatsu, Shiga 525-8577, Japan}
\email{takayama@se.ritsumei.ac.jp}
\date{March 6 in 2005}

\def\Coh#1#2{H_{\mm}^{#1}(#2)}

\newcommand{\AppTh}{Theorem~\ref{approxtheorem} }
\def\da{\downarrow}
\newcommand{\ua}{\uparrow}
\newcommand{\namedto}[1]{\buildrel\mbox{$#1$}\over\rightarrow}
\newcommand{\bdel}{\bar\partial}
\newcommand{\proj}{{\rm proj.}}

\maketitle

\opn\Hilb{Hilb}

\begin{abstract}
We give a generalization of Hochster's formula for local cohomologies of
square-free monomial ideals to monomial ideals, which are not
necessarily square-free.  Using this formula, we give combinatorial
characterizations of generalized Cohen-Macaulay monomial ideals.  We
also give other applications of the generalized Hochster's formula.
\end{abstract}

\section*{Introduction}

Let $K$ be a field and let $S = \poly{K}{X}{n}$ be a polynomial ring
with the standard grading.  For a graded ideal $I\subset S$ we set $R =
S/I$.  We denote by $x_i$ the image of $X_i$ in $R$ for $i=1,\ldots, n$
and set $\mm =(x_1,\ldots, x_n)$, the unique graded maximal ideal.  Also
$\Coh{i}{R}$ denotes the local cohomology module of $R$ with regard
$\mm$. A residue class ring $R$ is called a {\em generalized
Cohen-Macaulay} ring ({\em generalized CM} ring), 
or {\em FLC} (Finite Length Cohomology) ring,
if $\Coh{i}{R}$ has finite length  for $i\ne \dim R$.
In this case, we will call the ideal $I\subset S$ a 
{\em generalized CM ideal}.

As defining ideals of algebraic sets, we can find many examples of
generalized CM ideals such as homogeneous coordinate rings 
of non-singular projective varieties.
For monomial ideals, which are not directly
related to algebraic sets, the notions of generalized CM rings and 
Buchsbaum rings \cite{SV} coincide in the square-free case and 
the combinatorial characterization 
of generalied CM square-free monomial ideals (Stanley-Reisner ideals) 
has been given in terms Buchsbaum simplicial complexes \cite{Sch,St,SV}.  
However, as far as the author is concerned, the case of
non-square-free monomial ideals has not been studied very much, and the
aim of this paper is to give  combinatorial characterizations of
generalized CM monomial ideals, which are not always square-free.

We first give a generalization of Hochster's formula 
on local cohomologies for square-free monomial ideals \cite{Hoc}
to monomial ideals that are not necessarily square-free
(Theorem~\ref{hochster:theoremTheFormula}). 
From this formula, we can easily deduce several already known and 
probably new facts on vanishing degrees of local cohomologies. 
In particular, the vanishing degrees of generalized CM monomial 
ideals (Proposition~\ref{hochster:propFinLen}). 
This result allows us to deduce  combinatorial 
characterizations of generalized CM monomial ideals in terms of 
the exponents of variables in the monomial generators
(Theorem~\ref{finlen-main}, Corollary~\ref{2dim} and 
Theorem~\ref{finlen-3dim}).
On the other hand, thanks to the generalized Hochster's formula 
we can compare local cohomologies for $I$ and its radical $\sqrt{I}$
(Proposition~\ref{htt:compare} and~\ref{htt:genCM}), which,
together with the combinatorial characterization of generalized CM
property, suggests a method to construct generalized CM monomial ideals 
from Buchsbaum Stanley-Reisner ideals. Namely, by changing a square-free
generator $X_{i_1}\cdots X_{i_\ell}$ of a Buchsbaum Stanley-Reisner
ideal $J$ to a monomial $X_{i_1}^{a_{i_1}}\cdots X_{i_\ell}^{a_\ell}$
 $(a_j\in\NN, j=1,\ldots,\ell)$, we make a generator of a
generalized CM monomial ideal $I$ with $\sqrt{I}=J$, 
and the combinatorial characterization
shows a right choice of the exponents $a_j$.

One way of the construction is changing
all the occurrences of the variable $X_i$ in the minimal set of
generators
to $X_i^{a_i}$ with a fixed exponent $a_i$, $i=1,\ldots, n$ (Example~\ref{frobenius}).
In some specific case, we can show, using our combinatorial
characterization, that this is the only way of construction
 (Example~\ref{onlyFrobenius}).

For a finite set $S$ we denote by $\mid S\mid$ the cardinarity of 
$S$, and, for sets $A$ and $B$, $A\subset B$ means that 
$A$ is a subset of $B$, which may be equal to $A$.
\thanks{The author thanks J\"urgen Herzog for 
valuable discussions and detailed comments on the early version of 
the paper.}

\section{Local cohomologies of monomial ideals}

\subsection{Generalized Hochster's Formula}
In this subsection, we give a natural extension of Hochster's formula
on local cohomologies of Stanley-Reisiner ideals to 
monomial ideals. The proof goes along almost the same line 
as that for Stanley-Reisner ideals given, for example, 
in \cite{BH} chapter~5.3. But we will give 
a full detail for the readers' convenience.

Let $I\subset S$  be a monomial ideal,
which is not necessarily square-free.
Then we have 
\begin{equation*}
	\Coh{i}{R}\iso H^i(C^\bullet)
\end{equation*}
where $C^\bullet$ is the $\check{C}$ech complex defined as follows:
\begin{equation*}
 C^\bullet : 
0\To C^0 \To C^1 \To \cdots \To C^n \To 0,
\qquad C^t = \Dirsum_{1\leq i_1<\cdots <i_t \leq n}R_{x_{i_1}\cdots x_{i_t}}.
\end{equation*}
and the differential $C^t\To C^{t+1}$ of this complex is induced by
\begin{equation*}
(-1)^s nat : R_{x_{i_1}\cdots x_{i_t}} \To R_{x_{j_1}\cdots x_{j_{t+1}}}
\qquad \text{with } \{i_1,\ldots, i_t\} = \{j_1,\ldots, \hat{h}_s, \ldots, j_{t+1}\}
\end{equation*}
where $nat$ is the natural homomorphism to localized rings
and $R_{x_{i_1}\cdots x_{i_t}}$, for example, denotes localization by
$x_{i_1},\ldots, x_{i_t}$.

We can consider a $\ZZ^n$-grading to $\Coh{i}{R}$, $C^\bullet$
and $R_{x_{i_1}\cdots x_{i_t}}$ induced by the multi grading of 
$S$. See for example \cite{BH} for more detailed information about this complex.

Now we will consider the degree $a$ subcomplex $C^\bullet_a$ 
of $C^\bullet$ for any $a\in\ZZ^n$. Before that we will prepare the notation.
For a monomial ideal $I\subset S$, we denote by $G(I)$ the 
minimal set of monomial generators. Let $u = X_1^{a_1}\cdots X_n^{a_n}$ 
be a monimial with $a_i\geq 0$ for all $i$, 
then we define $\nu_j(u) = a_j$ for $j=1,\ldots, n$,
and $\supp(u)= \{i\;\vert\; a_i\ne 0\}$.
We set
$G_a = \{i \mid a_i < 0\}$ and $H_a = \{i \mid a_i > 0\}$
for $a\in\ZZ^n$.

\def\hochster#1{for all $u\in G(I)$ there exists $j\notin {#1}$
such that $\nu_j(u) > a_j \geq 0$}

\begin{Lemma}
\label{hochster:lemma1}
Let $x= x_{i_1}\cdots x_{i_r}$ with $i_1 < \cdots <i_r$ and 
set $F = \supp(x)$.
For all $a\in\ZZ^n$ we have $\dim_K (R_x)_a \leq 1$ and 
the following are equivalent
\begin{enumerate}
	\item [$(i)$] $(R_x)_a \iso K$ 
	\item [$(ii)$] 
	$F \supset G_a$ and \hochster{F}.
	\end{enumerate}
\end{Lemma}
Notice that the condition $a_i \geq 0$ 
in $(ii)$ is redundant because this  follows from the condition $F\supset G_a$.
But it is written for the readers' convenience.
\begin{proof} The proof of $\dim_K (R_x)_a \leq 1$ is  
verbatim the same as that of Lemma~5.3.6 (a) in \cite{BH}. 
Now we assume $(i)$, i.e.,
$(R_x)_a \ne 0$. This is equivalent to the condition that 
there exists a monomial $\sigma\in R$ and $\ell\in\NN$ such that 
\begin{enumerate}
\item [$(a)$] $x^m \sigma\ne 0$ for all $m\in\NN$, and 
\item [$(b)$] $\deg\displaystyle{\frac{\sigma}{x^\ell}} = a$,
\end{enumerate}
where $\deg$ denotes the multidegree.
We know from $(b)$ that we have $F \supset G_a$
because a negative degree $a_i (<0)$ in $a$ must 
come from the denominator of the fraction $\sigma/x^\ell$ and 
$F = \supp(x^\ell)$.
Now we know that $(a)$ is equivalent to the following 
condition:
for all $u\in G(I)$ and for all $m\in\NN$ we have
$u\not| (X_{i_1}^m\cdots X_{i_r}^m)(X_1^{b_1}\cdots X_n^{b_n})$
where we set $\sigma = x_1^{b_1}\cdots x_n^{b_n}$ with some integers $b_j\geq 0$,
$j=1,\ldots, n$. 
Namely,
for all $u\in G(I)$ there exists $i\notin F$
such that $\nu_i(u) > b_i$. Furthermore, 
we know from the condition $F \supset G_a$ 
that we have $a_i = b_i$ for $i\notin F$ since
by $(b)$ non-negative degrees in $a$ must come from $\sigma$.
Consequently we obtain $(ii)$.

Now we show the converse. Assume that we have $(ii)$. 
Set $\tau = \prod_{i\in H_a}x_i^{a_i}$   
and $\rho = \prod_{i\in G_a}x_i^{-a_i}$. 
Then since $F \supset G_a$ there exists $\ell \in \NN$ and 
a monomial $\sigma$ in $R$ 
such that 
\begin{equation}
\label{hochster:lemma1-eq1}
	x^\ell = \rho \sigma   
\end{equation}
Now we show that 
$\displaystyle{\frac{\sigma\tau}{x^\ell}}\ne 0$ 
in $R_x$.
$\displaystyle{\frac{\sigma\tau}{x^\ell}}\ne 0$ 
is equivalent to the condition
that 
$x^m(\sigma\tau)\ne 0$ 
for all $m\in\NN$. As in the above discussion,
this is equivalent to the condition
\begin{equation}
\label{hochster:lemma1-eq2}
	\text{for all } u \in G(I)\; \text{ there exists } i\notin F
\mbox{ such that } \nu_i(u) > b_i
\end{equation}
where we set 
$\sigma\tau = x_1^{b_1}\cdots x_n^{b_n}$ 
for some integers $b_j\geq 0$, $j=1,\ldots, n$.
But by (\ref{hochster:lemma1-eq1}) we have 
$i\notin\supp(\sigma)$ 
for
$i\notin F$, so that 
$b_i = \nu_i(\tau) = a_i(>0)$ 
(i.e., 
$i\in H_a$)
or $a_i = b_i = 0$ (i.e., $i\notin H_a\cup G_a$). Hence
we can replace ``$\nu_i(u)>b_i$'' in (\ref{hochster:lemma1-eq2})
by ``$\nu_i(u)>a_i \geq 0$'' and then 
(\ref{hochster:lemma1-eq2}) 
is assured by the assumption.
Thus we have 
$\displaystyle{\frac{\sigma\tau}{x^\ell}}\ne 0$ in $R_x$.
Therefore
\begin{equation*}
   \deg\displaystyle{\frac{\sigma\tau}{x^\ell}}
=  \deg\displaystyle{\frac{\sigma\tau}{\rho\sigma}}
=  \deg\prod_{i\in H_a\cup G_a}x_i^{a_i}
= \deg x^a = a
\end{equation*}
as requied.
\end{proof}

Let $a\in\ZZ^n$. By Lemma~\ref{hochster:lemma1} we see that 
$(C^i)_a$ has a basis
\begin{equation*}
\left\{
	b_F 
	: \begin{array}{l}
           F\supset G_a,\; \vert F\vert = i,\\
           \text{ and } \mbox{\hochster{F}}
          \end{array}
\right\}.
\end{equation*}
Restricting the differentation of $C^\bullet$ to the
$a$th graded piece, we obtain a complex $(C^\bullet)_a$
of finite dimensional $K$-vector spaces with differentation
$\partial : (C^i)_a\To (C^{i+1})_a$ given
by $\partial(b_F) = \sum(-1)^{\sigma(F, F')}b_{F'}$
where the sum is taken over all $F'$ such that 
$F' \supset F$ with $\vert F'\vert = i+1$
and  \hochster{F'}. Also we define $\sigma(F, F') =s$
if  $F'= \{j_0,\ldots, j_i\}$ and 
$F = \{j_0,\ldots, \hat{j}_s,\ldots, j_i\}$.
Then we describe the $a$th component of the local cohomology
in terms of this subcomplex:
$\Coh{i}{R}_a \iso H^i(C^\bullet) = H^i(C^\bullet_a)$.

Now we fix our notation on simplicial complex. A simplicial 
complex $\Delta$ on a finite set $[n] = \{1,\ldots, n\}$
is a collection of subsets of $[n]$ such that $F\in\Delta$
whenever $F\subset G$ for some $G\in \Delta$.
Notice that, for the convenience in the 
later discussions, we do not assume the condition
that $\{i\}\in \Delta$ for $i=1,\ldots, n$.
We define $\dim F = i$ if $\mid F\mid = i+1$
and $\dim \Delta = \max\{\dim F \mid F\in\Delta\}$,
which will be called the dimension of $\Delta$.
If we assume a 
linear order on $[n]$, say $1<2<\cdots<n$, then we will call
$\Delta$ {\em oriented}, and in this case we always denote
an element $F=\{i_1,\ldots, i_k\}\in\Delta$ 
with the orderd sequence $i_1 <\ldots < i_k$.
For a given oriented simplicial complex
of dimension $d-1$, 
we denote by ${\cal C}(\Delta)$ the augumented oriented chain complex 
of $\Delta$:
\begin{equation*}
{\cal C}(\Delta) :
0\To {\cal C}_{d-1}\overset{\partial}{\To}
     {\cal C}_{d-2}\To\cdots\To
     {\cal C}_0\overset{\partial}{\To}
     {\cal C}_{-1}\To 0
\end{equation*}
where 
\begin{equation*}
   {\cal C}_i = \Dirsum_{F\in\Delta, \dim F=i}\ZZ F
\qquad\text{and}\qquad
   \partial F = \sum_{j=0}^i(-1)^jF_j
\end{equation*}
for all $F\in \Delta$. Here we define
$F_j = \{i_0,\ldots, \hat{i}_j,\ldots, i_k\}$
for $F = \{i_0,\ldots, i_k\}$.
Now for an abelian group $G$, we define
the $i$th reduced simplicial homology $\tilde{H}_i(\Delta; G)$ 
of $\Delta$ to be the $i$th homology 
of the complex ${\cal C}(\Delta)\tensor G$ for 
all $i$. Also we define the 
$i$th reduced simplicial cohomology 
$\tilde{H}^i(\Delta; G)$ of $\Delta$
to be the $i$th cohomology
of the dual chain complex $\Hom_\ZZ({\cal C}(\Delta), G)$
for all $i$.
Notice that we have 
\begin{equation*}
\tilde{H}_{-1}(\Delta; G) 
 = \left\{
	\begin{array}{ll}
	   G  &  \text{if $\Delta = \{\emptyset\}$} \\
	   0  &  \text{otherwise}
	\end{array}
   \right.,
\end{equation*}
and if $\Delta = \emptyset$ then $\dim \Delta = -1$ and 
$\tilde{H}_{i}(\Delta; G)=0$ for all $i$.

Now we will establish an isomorphism between the complex
$(C^\bullet)_a$, $a\in\ZZ^n$, and a dual chain complex.
For any $a\in\ZZ^n$, we define a simplicial complex
\begin{equation*}
\Delta_a = 
\left\{ F- G_a
	\;\vert\;
	\begin{array}{l}
	 F \supset G_a, \mbox{ and }\\
        \mbox{\hochster{F}}
	\end{array}
\right\}.
\end{equation*}
Notice that we may have $\Delta_a=\emptyset$ for some $a\in\ZZ^n$.

\begin{Lemma}
\label{hochster:lemma2}
For all $a\in\ZZ^n$ there exists an isomorphism
of complexes
\begin{equation*}
  \alpha^\bullet 
:  (C^\bullet)_a \To 
   \Hom_\ZZ({\cal C}(\Delta_a)[-j-1], K)
\qquad j = \vert G_a\vert
\end{equation*}
where ${\cal C}(\Delta_a)[-j-1]$ means shifting the degree 
of ${\cal C}(\Delta_a)$ by $-j-1$.
\end{Lemma}
\begin{proof}
The assignment $F\mapsto F - G_a$ induces an isomorphism
$\alpha^\bullet : (C^\bullet)_a \To \Hom_\ZZ({\cal C}(\Delta_a)[-j-1], K)$ 
of $K$-vector spaces such that $b_F \mapsto \varphi_{F-G_a}$,
where
\begin{equation*}
   \varphi_{F'}(F^")
= \left\{
	\begin{array}{ll}
	    1 & \mbox{if $F' = F^"$} \\
            0 & \mbox{otherwise.}
	\end{array}
  \right. 
\end{equation*}
That this is a homomorphism of complexes can be checked 
in a straightforward way.
\end{proof}

Now we come to the  main theorem in this section.

\begin{Theorem}
\label{hochster:theoremTheFormula}
Let $I\subset S = \poly{K}{X}{n}$ be a monomial ideal. Then
the multigraded Hilbert series of the local cohomology modules of $R = S/I$ 
with respect to the $\ZZ^n$-grading is given by
\begin{equation*}
\Hilb(\Coh{i}{R}, {\bf t})
=\sum_{F\in\Delta}
  \sum
   \dim_K\tilde{H}_{i-\vert F\vert -1}(\Delta_a; K) {\bf t}^a
\end{equation*}
where 
${\bf t} = t_1\cdots t_n$,
the second sum runs over 
$a\in\ZZ^n$ such that 
$G_a = F$ and $a_j\leq \rho_j-1$, $j=1,\ldots,n$,
with 
$\rho_j = \max\{\nu_j(u)\;\vert\; u\in G(I)\}$
for $j=1,\ldots, n$, and $\Delta$ is the 
simplicial complex corresponding to the Stanley-Reisner
ideal $\sqrt{I}$.
\end{Theorem}

\begin{proof}
By Lemma~\ref{hochster:lemma2} and universal coefficient theorem
for simplicial (co)homology, we have 
\begin{eqnarray*}
\Hilb(\Coh{i}{R},{\bf t})
& = & \sum_{a\in\ZZ^n}\dim_K\Coh{i}{R}_a {\bf t}^a 
  =  \sum_{a\in\ZZ^n}\dim_K H^i(C^\bullet_a) {\bf t}^a \\
& = & \sum_{a\in\ZZ^n}\dim_K\tilde{H}_{i-\vert G_a\vert -1}(\Delta_a;K) 
                                         {\bf t}^a. \\
\end{eqnarray*}
It is clear from the definition that $\Delta_a= \emptyset$
if for all $j\notin G_a$ we have $a_j \geq \rho_j$.
Moreover for all $a\in\ZZ^n$ with 
$a_j\geq \rho_j$ for at least one index $j\notin G_a$
we have $\dim_K\tilde{H}_{i-\vert G_a\vert -1}(\Delta_a;K) =0$.
To prove this fact we can assume without loss of generality that 
$a_1\geq \rho_1$ and that $\Delta_a\ne\emptyset$. Then 
we have $1\notin G_a$, and, for all $\sigma = (L - G_a)\in \Delta_a$
with $L \supset G_a$ and 
$1\notin \sigma$, we have  $\sigma\cup \{1\}\in \Delta_a$.
In fact, since we have $\nu_1(u)\leq a_1$ for all $u\in G(I)$
 the existence of 
$k\notin L$ with $\nu_k(u)>a_k$ implies $k\notin L\cup\{1\}$.
Consequently we know that $\Delta_a$ is a cone by the vertex $\{1\}$ 
so that, as is well known, we have 
$\tilde{H}_{i-\vert G_a\vert -1}(\Delta_a;K) =0$
for all $i$ as required.
Thus we obtain
\begin{eqnarray*}
\Hilb(\Coh{i}{R},{\bf t})
& = & \sum_{{\tiny
		\begin{array}{c}
	              a\in\ZZ^n\\
		   a_j\leq \rho_j-1\\
	           j=1,\ldots, n
		\end{array}
             }}
        \dim_K\tilde{H}_{i-\vert G_a\vert -1}(\Delta_a;K) 
                                         {\bf t}^a. \\
\end{eqnarray*}
Now if $\Delta_a\ne \emptyset$, we must have $(G_a - G_a =)\;
\emptyset\in \Delta_a$, i.e., 
\hochster{G_a}, and this implies that 
$G_a \not\supset \supp(u)$ for all $u\in G(I)$,
namely $G_a$ is not 
a non-face of $\Delta$, i.e., $G_a\in \Delta$.
Thus we finally obtain the required formula.
\end{proof}

The original Hochster's formula is a special case of 
Theorem~\ref{hochster:theoremTheFormula}. 

\begin{Corollary}[Hochster]
Let $\Delta$ be a simplicial complex and 
let $K[\Delta]$ be the Stanley-Reisner ring 
corresponding to $\Delta$. Then we have
\begin{equation*}
\Hilb(\Coh{i}{K[\Delta]}, {\bf t})
=\sum_{F\in\Delta}
   \dim_K\tilde{H}_{i-\vert F\vert -1}(\lk_\Delta F; K) 
	\prod_{j\in F}
	   \displaystyle{\frac{t_j^{-1}}{1- t_j^{-1}}},
\end{equation*}
where $\lk_\Delta F = \{G \vert F\cup G\in\Delta, F\cap G=\emptyset\}$.
\end{Corollary}
\begin{proof}
By Theorem~\ref{hochster:theoremTheFormula} we have 
\begin{equation*}
\Hilb(\Coh{i}{R}, {\bf t})
=\sum_{F\in\Delta}
 \sum_{{\tiny \begin{array}{c}
	a\in\ZZ^n_{-} \\
        G_a = F
	\end{array}}}
   \dim_K\tilde{H}_{i-\vert F\vert -1}(\Delta_a; K) {\bf t}^a
\end{equation*}
where $\ZZ^n_{-} = \{a\in\ZZ^n \vert a_j\leq 0\text{ for }j=1,\ldots,n\}$
and 
\begin{eqnarray*}
\Delta_a &=& 
\left\{ F- G_a
	\;\vert\;
	\begin{array}{l}
	 F \supset G_a, \mbox{ and for all $u\in G(I)$ there exists $j\notin F$}\\
        \mbox{
           such that $j\in \supp(u)$ and $j\notin H_a\cup G_a$
             }
	\end{array}
\right\}.\\
&=& 
\left\{ F- G_a
	\;\vert\;
	 F \supset G_a, \text{ and for all $u\in G(I)$ we have
          $H_a\cup F\not\supset \supp(u)$}
\right\}.\\
&=& 
\left\{ L
	\;\vert\;
         L\cap G_a=\emptyset, 
	 L\cup G_a\cup H_a\in\Delta
\right\} = \lk_{\st H_a}G_a.
\end{eqnarray*}
Then the rest of the proof is exactly as in Theorem~5.3.8 \cite{BH}.
\end{proof}

\subsection{Vanishing degrees of local cohomolgies}

In this subsection, we give some easy consequences  of 
Theorem~\ref{hochster:theoremTheFormula}. 
We define $a_i(R) = \max\{j \vert \Coh{i}{R}_j\ne 0\}$
if $\Coh{i}{R}\ne 0$ and $a_i(R) = -\infty$ if $\Coh{i}{R}=0$.
Similarly, we define 
and $b_i(R) = \inf\{j \vert \Coh{i}{R}_j\ne 0\}$
if $\Coh{i}{R}\ne 0$ and  $b_i(R) = +\infty$ if 
$\Coh{i}{R} =0$.

Recall that 
$\rho_j = \max\{\nu_j(u)\;\vert\; u\in G(I)\}$
for $j=1,\ldots, n$.

\begin{Corollary}
\label{hochster:corollary1}
Let $I\subset S = \poly{K}{X}{n}$ be a monomial ideal.
Then $a_i(R) \leq \sum_{j=1}^n \rho_j - n$ for all $i$.
\end{Corollary}
\begin{proof}
By Theorem~\ref{hochster:theoremTheFormula},
the terms in 
$\Hilb(\Coh{i}{R}, {\bf t})$ 
with the highest total degree are 
at most $\dim_K\tilde{H}_{i-\vert F\vert -1}(\Delta_a;K){\bf t}^a$ 
with $a_j = \rho_j-1$ for $j=1,\ldots, n$. Thus the total degree is 
at most $\sum_{j}\rho_j - n$.
\end{proof}

From Corollary~\ref{hochster:corollary1}, we can recover the 
following well known result.

\begin{Corollary}
Let $I\subset S$ be a generalized CM Stanley-Reisner ideal. Then 
$a_i(R) \leq 0$ for all $i$.
\end{Corollary}
\begin{proof}
If $I$ is square-free, then $\rho_j\leq 1$ for $j=1,\ldots, n$.
\end{proof}

For a Stanley-Reisner generalized CM ideal $I\subset S$
with $\dim R =d$,
it is well known that it is Buchsbaum and $b_i(R)\geq 0$ for all $i(\ne d)$. 
The following theorem extends this result to monomial ideals in general.

\begin{Proposition}
\label{hochster:propFinLen}
Let $I\subset S=\poly{K}{X}{n}$ be a monomial ideal.
Then following are equivalent:
\begin{enumerate}
\item [$(i)$]  $\length(\Coh{i}{S/I})<\infty$
\item [$(ii)$] $\Coh{i}{S/I}_a=0$ for all $a\in\ZZ^n$
with $G_a\ne\emptyset$, in particular $b_i(S/I)\geq 0$
\item [$(iii)$] $\tilde{H}_{i-\vert G_a\vert -1}(\Delta_a; K)=0$ 
for all $a\in\ZZ^n$ with $a_j\leq \rho_j-1$
      $(j=1,\ldots, n)$ and $\emptyset\ne
G_a\in\Delta$.
\end{enumerate}
\end{Proposition}
\begin{proof}
The equivalence of $(ii)$ and $(iii)$ are immediate from
Theorem~\ref{hochster:theoremTheFormula}. We will prove the 
equivalence of $(i)$ and $(iii)$.
Assume that $\length(\Coh{i}{S/I})<\infty$.
Assume also that there exists $a\in\ZZ^n$ such that 
$a_j\leq \rho_j-1$ $(j=1,\ldots,n)$,
$\emptyset\ne G_a\in\Delta$ and $\tilde{H}_{i-\vert G_a\vert -1}(\Delta_a;K)\ne 0$.
Now observe that by the definition of $\Delta_a$, the condition
is independent of the values  $a_j$ for $j\in G_a$.
This means that the total degree $j = \sum_{k=1}^n a_k$ can be 
any negative integer so that $\Coh{i}{R}$ is not of finite length,
which contradicts the assumption.
Thus we must have $\tilde{H}_{i-\vert G_a\vert -1}(\Delta_a; K)=0$
for all such $a\in\ZZ^n$.
The converse implication is straightforward.
\end{proof}

\begin{Corollary}
\label{hochster:theoremGenCMCase}
Let $I\subset S=\poly{K}{X}{n}$ be a generalized CM
monomial ideal with $\dim R = d(>0)$.
Then $b_i(R) \geq 0$ for all $i (\ne d)$.
\end{Corollary}

For a generalized CM ring $R$,
there exists an integer $k\in\ZZ$,
$k\geq 1$, such that $\mm^k\Coh{i}{R} = 0$ for $i\ne \dim R$. If this
condition holds, we will also call $R$, or $I\subset S$, {\em
$k$-Buchsbaum}. An ideal $I$ is generelized CM if and only if it is
$k$-Buchsbaum for some $k$. 
If $I$ is $k$-Buchsbaum but not
$(k-1)$-Buchsbaum, then we will call $I$ {\em strict $k$-Buchsbaum}.

\begin{Proposition}
\label{hochster:theoremTheBound}
Let $I\subset S=\poly{K}{X}{n}$ be a generalized CM
monomial ideal. Then $R=S/I$ is 
$\left(\sum_{j=1}^n \rho_j -n+1\right)$-Buchsbaum.
\end{Proposition}
\begin{proof}
$R$ is $\max\{a_i(R) - b_i(R) + 1\;\vert\; i\ne d\}$-Buchsbaum.
Then the required result follows immediately from
Corollary~\ref{hochster:corollary1} and 
Corollary~\ref{hochster:theoremGenCMCase}.
\end{proof}

From Proposition~\ref{hochster:theoremTheBound}, we immediately
know that a Stanley-Reisner ideal is 1-Buchsbaum if it is 
generalized CM, which is a weaker version of 
the well-known result that a generalized CM Stanley-Reisner 
ideal is Buchsbaum.


The bound of $k$-Buchsbaumness given in
Proposition~\ref{hochster:theoremTheBound} is best possible.
In fact, we can construct strict $(\sum_{j=1}^n \rho_j -n+1)$-Buchsbaum
ideals  as in the following example.

\begin{Example}
{\em
\label{frobenius}
Let $I\subset S$ be a Stanley-Reisner Buchsbaum ideal.
Notice that such ideals can be constructed with the method
presented in \cite{BjoHib} and $\Coh{i}{S/I}$ ($i\ne \dim R$)
is a $K$-vector space for $i\ne \dim R$.

Now consider a $K$-homomorphism 
\begin{equation*}
  \varphi: S \To S, \qquad  X_i \longmapsto X_i^{a_i}\; (i=1,\ldots, n)
\end{equation*}
where $a = (a_1,\ldots, a_n) \in\ZZ^n$ with $a_i\geq 1$ for $i=1,\ldots, n$.
We define $\varphi(M) = M\tensor_S \hbox{}^\varphi S$ for a $S$-module,
where a left-right $S$-module $\hbox{}^\varphi S$ is equal to $S$ as a set,
it is a right $S$-module in the ordinary sense and 
its left $S$-module structure is determined by $\varphi$.
Then we have 
\begin{enumerate}
\item $\varphi(S/I) = S/\varphi(I)S$,
\item $\varphi$ is an exact functor.
\end{enumerate}
Thus, for $i\ne \dim R$, 
we have $\Coh{i}{S/\varphi(I)S}\iso \varphi(\Coh{i}{S/I})$ and 
since $\Coh{i}{S/I}$ is a direct sum of $S/\mm$, 
$\Coh{i}{S/\varphi(I)S}$ is a direct sum of 
$S/(X_1^{a_1},\ldots, X_n^{a_n})$. Then we know that 
$\mm^k\Coh{i}{S/I} = 0$ but  $\mm^{k-1}\Coh{i}{S/I}\ne 0$ 
with $k = \sum_{j=1}^n \rho_j - n +1 = \sum_{j=1}^n a_j - n +1$.
}
\end{Example}

\begin{Remark}
{\em 
Bresinsky and Hoa gave a bound for $k$-Buchsbaumness
for ideals generated by monomials and binomials 
$($Theorem~4.5 \cite{BreHoa}$)$.
For monomial ideals, our bound is stronger than that of 
Bresinsky and Hoa. 
Also, according to K.\ Yanagawa, 
Proposition~\ref{hochster:theoremTheBound}
can also be deduced from his theory of square-free modules \cite{Y}.
}
\end{Remark}

Recall that Castelnuovo-Mumford regularity of the ring $R$ 
is defined by 
\begin{equation*}
   \reg(R) = \max\{i+j \vert \Coh{i}{R}_j\ne 0\}.
\end{equation*}
Let $r = \reg(R)$. Then we have $\Coh{i}{R}_j =0$ 
for $j > r -i$. 
Then we have

\begin{Corollary}
Let $I\subset S$ be a generalized CM monomial ideal
with $d =\dim R$ and  $r = \reg(R)$. 
Then $\Coh{i}{R} =0$ for $r+1 \leq i<d$.
In particular, if $I$ has $q$-linear resolution, 
we have  $\Coh{i}{R} =0$ for $q \leq i<d$.
\end{Corollary}
\begin{proof}
First part is clear from Corollary~\ref{hochster:theoremGenCMCase}.
If $R$ has $q$-linear resolution, we have $\reg(R) = q-1$.
Thus the last part also follows immediately.
\end{proof}

\section{Generalized Cohen-Macaulay monomial ideals}

\subsection{FLC property}
In this subsection, we give a combinatorial characterization
of FLC (finite length cohomology) property for monomial ideals,
as an application of Theorem~\ref{hochster:theoremTheFormula}.
We prepare some notations. 
Let $I\subset S = \poly{K}{X}{n}$ be a monomial ideal.  If $X_i^a\in
G(I)$ for some $1\leq i\leq n$ and $a\in\NN$, we easily know that $a$
must be $\rho_i=\max\{\nu_i(u)\;\vert\; u\in G(I)\}$.  Then, by changing
the name of the variables if necessary, we can write without loss of
generality that $G(I) = \{X_{m+1}^{\rho_{m+1}}, \ldots, X_n^{\rho_n}\}
\cup G_0(I)$ with $m\leq n$, where $G_0(I) = \{u\in G(I) \;\vert\;
\vert\supp(u)\vert\geq 2\}$.
We denote by $\Delta$ the simplicial complex corresponding to 
a square-free monomial ideal $\sqrt{I}$, which is a complex
over the vertex set $[m] = \{1,\ldots, m\}$.
We regard $\ZZ^n$ as a partially ordered set
by defining $a\leq b$, $a, b\in\ZZ^n$,  to be 
$a_i\leq b_i$ for $i=1,\ldots, n$.
We denote $\rho-1 = (\rho_1-1,\ldots,
\rho_n-1)\in\ZZ^n$. 
For $a\in\ZZ^n$ and a monomial $u\in S$, we define
$L(a,u) = \{ i\in [n]\;\vert\; \nu_i(u)> a_i\}$.
Also for $a\in\ZZ^n$ with $a\leq \rho -1$ and $\sigma\subset [n]$,
we define $a(\sigma)\in\ZZ^n$ as follows:
\begin{equation*}
    a(\sigma)_i = 
\left\{
 \begin{array}{ll}
  a_i    & \mbox{if $i\notin\sigma$} \\
  \rho_i & \mbox{if $i\in\sigma$}
 \end{array}
\right.
\end{equation*}
We abbreviate $a(\{j\})$ as $a(j)$ for $j\in\ZZ$.


Now we prove
\begin{Theorem}
\label{finlen-main}
Let $I \subset S = \poly{K}{X}{n}$ be a monomial ideal.
If $\length(\Coh{i}{S/I}) < \infty$ $(i>0)$ then,
for all $(i-1)$-face $\sigma\in\Delta$
and for all $a\in\ZZ^n$
such that 
\begin{enumerate}
\item [$(a)$] $0\leq a\leq \rho -1$, 
\item [$(b)$] $L(a(\sigma), u)\ne \emptyset$ for all $u\in G_0(I)$, and 
\item [$(c)$] $a(\sigma)$ is maximal with the properties $(a)$ and $(b)$
for a fixed $\sigma$, 
\end{enumerate}
we have the following:
there exists  $\ell\in [m]\backslash \sigma$ such that
\begin{enumerate}
\item [(i)]  $a_\ell = \rho_\ell -1$, and 
\item [(ii)] for all $u\in G_0(I)$ with $\nu_\ell(u) =\rho_\ell$
we have $L(a(\sigma\cup\{\ell\}), u)\ne \emptyset$.
\end{enumerate}
\end{Theorem}
\begin{proof}Assume that $\length(\Coh{i}{S/I})<\infty$. Then,
by Proposition~\ref{hochster:propFinLen}, 
we have $\tilde{H}_{-1}(\Delta_a;K)=0$ 
for all $a\in\ZZ^n$ with $a\leq \rho-1$, $G_a\in\Delta$
and $\vert G_a\vert =i(>0)$. This implies $\Delta_a\ne\emptyset$.
Now for such $a\in\ZZ^n$ we set $\sigma = G_a$.  Notice that $\sigma
\subset [m]$ since $\sigma\in\Delta$. 
We also notice that, as far as the complex 
$\Delta_a$ is concerned, the values $a_i$ for $i\in\sigma=G_a$ are irrelevant. Thus
we will change the values $a_i$ $(i\in\sigma)$ and 
assume that $0\leq a\leq \rho-1$. Notice that 
$\Delta_a$, for the new $a$, is the same as before we change 
the values $a_i$ for $i\in\sigma$.

We know that $\Delta_a = \emptyset$ is equivalent to the
condition that there exists $u\in G(I)$ such that 
$L(a(\sigma),u) =\emptyset$. Now we assume 
that $L(a(\sigma),u) \ne \emptyset$ for all $u\in G(I)$, namely
$\Delta_a\ne\emptyset$.

The condition $\Delta_a\ne \{\emptyset\}$ 
is equivalent to the
condition that $\{\ell\}\in \Delta_a$ for some $\ell\in [n]$, i.e.,
there exists $\ell\in [n]\backslash\sigma$ such that for all $u\in G(I)$
we can find  $k\in [n]\backslash(\sigma\cup\{\ell\})$ satisfying
$\nu_k(u)>a_k$. 
Namely,
\begin{equation}
\label{condition}
\mbox{there exists $\ell\in[n]\backslash\sigma$ such that }
L(a(\sigma\cup\{\ell\}), u)\ne \emptyset\; \mbox{for all }u\in G(I).
\end{equation}
Under this condition 
we have,  for any $b\in\ZZ^n$ with $0\leq b\leq \rho-1$ and 
$b(\sigma)\leq a(\sigma)$,
$\L(b(\sigma\cup\{\ell\}), u)\supseteq 
 L(a(\sigma\cup\{\ell\}),u) \ne \emptyset$ for all $u\in G(I)$.
Thus we can assume that $a(\sigma)$ is maximal 
satisfying the condition that $L(a(\sigma),u) \ne \emptyset$ for 
all $u\in G(I)$ and $0\leq a\leq \rho-1$. Also,
since $\sigma \subset [m]$ and $a\leq \rho-1$,
we have  $L(a(\sigma),X_i^{\rho_i})=\{i\}\ne\emptyset$
for all $m+1\leq i\leq n$. Hence we can replace '$u\in G(I)$' by
'$u\in G_0(I)$' in the maximality condition for $a(\sigma)$.
Now we have only to show that the condition (\ref{condition})
is equivalent to $(i)$ and $(ii)$ in the statement.

Since we have $L(a(\sigma\cup\{\ell\}), X_j^{\rho_j})=\emptyset$
for all $m+1\leq j\leq n$, we can only  find the index $\ell$ 
as in $(\ref{condition})$ in $[m]\backslash\sigma$.  
Now for  $\ell\in [m]\backslash\sigma$, the
existence of $k\in [n]\backslash(\sigma\cup\{\ell\})$ satisfying
$\nu_k(u)>a_k$ is always assured for every $u\in
\{X_{m+1}^{\rho_{m+1}},\ldots, X_n^{\rho_n}\}$.  Moreover, if
$\ell\notin\supp(u)$ for $u\in G_0(I)$, then $L(a(\sigma\cup\{\ell\}),
u) = L(a(\sigma),u)$ and this is $\ne \emptyset$ since 
$\Delta_a\ne\emptyset$.  Thus $(\ref{condition})$ is equivalent
to the existence of $\ell\in [m]\backslash\sigma$ such that 
\begin{equation} 
\label{condition2}
L(a(\sigma\cup\{\ell\}), u) \ne \emptyset\mbox{ for }u\in G_0(I)
\mbox{ with }\ell\in\supp(u).
\end{equation}

Assume that $a_\ell < \rho_\ell -1$ and set 
$e\in\ZZ^n$ as $e_i=0$ for $i\ne \ell$ and $e_\ell = \rho_\ell-1 - a_\ell$.
Then by the maximality of 
$a(\sigma)$ there exists $u\in G_0(u)$ such that 
$\emptyset = L(a(\sigma) + e, u) \supset L(a(\sigma\cup{\ell}),u)$, 
which contradicts the condition $(\ref{condition})$. Thus we must 
have $a_\ell = \rho_\ell -1$ for $\ell$ as in $(\ref{condition})$.
If $u\in G_0(I)$ is such that $0< \nu_\ell(u)< \rho_\ell$, 
then $L(a(\sigma)\cup\{\ell\},u) = L(a(\sigma),u) \ne \emptyset$
by assumption on $a(\sigma)$. Thus we can replace
'$u\in G_0(I)$ with $\ell\in\supp(u)$' in the condition
(\ref{condition2}) by '$u\in G_0(I)$ with $\nu_\ell(u)=\rho_\ell$'.
Consequently we know that $(\ref{condition2})$ is equivalent to
$(i)$ and $(ii)$.
\end{proof}

From Theorem~\ref{finlen-main}, we can recover a weaker version of 
the well-known result as follows.

\begin{Corollary}
If $I\subset S$ is a generalized CM Stanley-Reisner ideal, i.e., Buchsbaum ideal,
then $\Delta$ is pure, namely, every facet has the same dimension.
\end{Corollary}
\begin{proof}
Let $I\subset S$ be a generalized CM Stanley-Reisner ideal. Then
by Theorem~\ref{finlen-main} we know that 
for every $0< i<\dim S/I$ and for every $(i-1)$-face $\sigma\in\Delta$
there exists $\ell$ with $1 \leq \ell \leq n$ such that
$\sigma\cup \{\ell\}\in \Delta$. From this we immediately know that 
$\Delta$ is pure.
\end{proof}

\subsection{generalized CM monomial ideals of $\dim\leq 3$}
If $\dim R \leq 1$, $I$ is always (generalized) CM. For $dim R=2, 3$, we 
can give combinatorial characterizations of generalized CM monomial
ideals
as follows.
First we give the $\dim 2$ case.

\begin{Corollary}
\label{2dim}
A monomial ideal $I\subset S$
is generalized CM with $\dim S/I = 2$ 
if and only if 
\begin{enumerate}
\item [$(i)$] $\dim \Delta=1$, and 
\item [$(ii)$]
for all $j\in [m]$ and for all $a\in\ZZ^n$ such that 
$(a)$ $0\leq a\leq \rho -1$, 
$(b)$ $L(a(j), u)\ne \emptyset$ for all $u\in G_0(I)$, and 
$(c)$ $a(j)$ is maximal with the properties $(a)$ and $(b)$
for a fixed $j$, 
we have the following:
there exists  $\ell\in [m]\backslash \{j\}$ such that
\begin{enumerate}
\item [(i)]  $a_\ell = \rho_\ell -1$, and 
\item [(ii)] for all $u\in G_0(I)$ with $\nu_\ell(u) =\rho_\ell$
we have $L(a(\{j,\ell\}), u)\ne \emptyset$.
\end{enumerate}
\end{enumerate}
\end{Corollary}
\begin{proof}
As is well-known, 
$\dim S/I=\dim S/\sqrt{I}= 2$ if and only if $\dim \Delta=1$.
Then $S/I$ is generalized CM if and only if 
$\length(\Coh{i}{S/I})<\infty$ for $i=0, 1$.
$\Coh{0}{S/I}$ is always of finite length and 
$\Coh{1}{S/I}$ is of finite length if and only if 
$\tilde{H}_{-\vert G_a\vert}(\Delta_a; K)=0$ 
for all $a\in\ZZ^n$ with $a\leq \rho-1$ and $\emptyset\ne
G_a\in\Delta$
by Proposition~\ref{hochster:propFinLen}. 
If $\vert G_a\vert \geq 2$,
we always have $\tilde{H}_{-\vert G_a\vert}(\Delta_a; K)=0$.
Now let $G_a = \{j\}$. Since $G_a\in\Delta$ we must 
have $j\in [m]$. Also $\tilde{H}_{-1}(\Delta_a; K)=0$
if and only if $\Delta_a\ne\emptyset$,
which is equivalent to the condition in the statement 
by the proof of Theorem~\ref{finlen-main}.
\end{proof}

In $\dim 3$ case, we need to give a combinatorial criterion
for connectedness of simplicial complexes. Notice that 
a simplicial complex $\Delta$ over the vertex set $[m]$ is 
{\em not} connected if and only if there exists non-empty 
disjoint subsets
$P, Q\subset [m]$ such that $P\cup Q = [m]$
and for all $p\in P$ and all $q\in Q$ there is no 1-face $\{p,q\}
\in\Delta$. 

\begin{Lemma}
\label{vset}
Assume that $G_a = \{j\}$, $j\in [m]$, and 
$\Delta_a\ne \emptyset$. Then
the set of vertices of $\Delta_a$ is 
$\{\ell\in [m]\backslash\{j\} \;\vert\;
\mbox{for all } u\in G_0(I)\mbox{ we have }
L(a(\{\ell, j\}), u) \ne \emptyset\}$.
\end{Lemma}
\begin{proof}
The $0$th skeleton of $\Delta_a$ is 
\begin{eqnarray*}
\lefteqn{
\{\{\ell\} \;\vert\; \ell\ne j,
\mbox{ for all } u\in G(I)\mbox{ there exists }k\notin \{\ell, j\}
\mbox{ such that }\nu_k(u)> a_k\}}\\
&=& \{\{\ell\} \;\vert\;  \ell\in [m]\backslash\{j\},
\mbox{ for all } u\in G_0(I)\mbox{ there exists }k\notin \{\ell, j\}
\mbox{ such that }\nu_k(u)> a_k\}\\
&=& \{\{\ell\} \;\vert\;  \ell\in [m]\backslash\{j\},
\mbox{ for all } u\in G_0(I)\mbox{ we have }
L(a(\{\ell, j\}), u) \ne \emptyset\},
\end{eqnarray*}
where the first equation is because if $m+1\leq \ell \leq n$
there is no $k\notin \{\ell, j\}$ such that $\nu_k(X_\ell^{\rho_\ell})
> a_k$ and if $\ell\in [m]$ we always have the index $k\notin\{\ell,j\}$
such that $\nu_k(X_i^{\rho_i})>a_i$ for $i=m+1,\ldots, n$, which is 
actually $k=i$.
\end{proof}

Now we show a combinatorial characterization of 
$\dim 3$ generalized CM monomial ideals.

\begin{Theorem}
\label{finlen-3dim}
A monomial ideal $I\subset S$ 
is  generalized CM with $\dim S/I =3$ if and only if 
\begin{enumerate}
\item [$(i)$] $\dim \Delta = 2$, and 
\item [$(ii)$] for all $j\in [m]$ 
and for all $a\in\ZZ^n$ such that 
$(a)$ $0\leq a \leq \rho -1$,
$(b)$ $L(a(j), u) \ne \emptyset$ for all $u\in G_0(I)$, and 
$(c)$ $a(j)$ is maximal with the properties $(a)$ and $(b)$
for a fixed $j$,
we have 
      \begin{enumerate}
       \item [$(1)$] there exists $\ell\in [m]\backslash\{j\}$ such that 
	     \begin{enumerate}
	      \item  [1.] $a_\ell = \rho_\ell -1$, and 
	      \item  [2.] for all $u\in G_0(I)$ with $\nu_\ell(u)=\rho_\ell$ we
		    have $L(a(\{j, \ell\}),u)\ne \emptyset$
	     \end{enumerate}
       \item [$(2)$] there are no non-empty disjoint subsets $P, Q\subset [m]$
	     satisfying the following property:
       \begin{enumerate}
	\item [1.] $P\cup Q = [m]\backslash\{j\}- L_a$\\
	      where $L_a = \{\ell \;\vert\; L(a(j), u)=\{\ell\} 
                  \mbox{ for some}\ u\in G_0(U) \}$, and 
	\item [2.] for all $x\in P$ and all $y\in Q$
	      there exists $u\in G_0(I)$ such that 
	      $L(a(j), u) = \{x,y\}$
       \end{enumerate}
      \end{enumerate}
\item [$(iii)$]
for all 1-face  $\sigma = \{i,j\}\in\Delta$
and for all $a\in\ZZ^n$
such that 
$(a)$ $0\leq a\leq \rho -1$, 
$(b)$ $L(a(\{i,j\}), u)\ne \emptyset$ for all $u\in G_0(I)$, and 
$(c)$ $a(\{i,j\})$ is maximal with the properties $(a)$ and $(b)$
for a fixed $\{i,j\}$, 
we have the following:
there exists  $\ell\in [m]\backslash \{i,j\}$ such that
\begin{enumerate}
\item [1.]  $a_\ell = \rho_\ell -1$, and 
\item [2.] for all $u\in G_0(I)$ with $\nu_\ell(u) =\rho_\ell$
we have $L(a(\{i,j,\ell\}), u)\ne \emptyset$.
\end{enumerate}
\end{enumerate}
\end{Theorem}
\begin{proof}
$\dim S/I =\dim S/\sqrt{I}=3$  if and only if $\dim\Delta=2$.
Now assume that $\dim S/I =3$. Then
$S/I$ is generalized CM if and only if 
$\length(\Coh{1}{S/I})<\infty$ and $\length(\Coh{2}{S/I})<\infty$,
which is equivalent to 
\begin{equation}
\label{1} \tilde{H}_{-\vert G_a\vert}(\Delta_a; K)=0 
\end{equation}
and 
\begin{equation}
\label{2} \tilde{H}_{1-\vert G_a\vert}(\Delta_a; K)=0 
\end{equation}
for all $a\in\ZZ^n$ with $a\leq \rho-1$ and $\emptyset\ne
G_a\in\Delta$ 
by Proposition~\ref{hochster:propFinLen}.
The condition $(\ref{1})$ 
is equivalent to $(ii)(1)$ by Corollary~\ref{2dim} 
and the condition $(\ref{2})$ 
is equivanent to 

\begin{equation}
\label{3} \tilde{H}_{0}(\Delta_a; K)=0
\mbox{ for all $a\in\ZZ^n$ with $a\leq \rho-1$ and $G_a = \{j\} \in [m]$},
\end{equation}
and 
\begin{equation}
\label{4} \tilde{H}_{-1}(\Delta_a; K)=0
\mbox{ for all $a\in\ZZ^n$ with $a\leq \rho-1$ and $G_a = \{i,j\} \in [m]$}
\end{equation}
%
since $\tilde{H}_{-k}(\Delta_a; K)=0$ for $k\geq 2$.
The condition $(\ref{3})$ 
exactly 
means the connectedness of the simplicial complex
$\Delta_a$. Let ${\cal V}_a$ be the set of vertices 
of $\Delta_a$.
By what we noticed just before Lemma~\ref{vset}, this is
equivalent to the condition that there exist
disjoint no non-empty 
subsets $P, Q\subset {\cal V}_a$ such that 
$P\cup Q = {\cal V}_a$
and for all $x\in P$ and all $y\in Q$ we have $\{x,y\}\notin
\Delta_a$. 
By Lemma~\ref{vset} we have 
${\cal V}_a = 
[m]\backslash\{j\} - L_a$ where
$L_a= \{\ell \;\vert\; \ell\ne j,\; L(a(\{j, \ell\}),u)=\emptyset
\mbox{ for some }u\in G_0(I)\}$. 
Since we have $L(a(j),u)\ne\emptyset$, $L(a(\{j,\ell\}),u)=\emptyset$
implies $L(a(j),u)=\{\ell\}$. Thus $L_a= \{\ell \;\vert\;
L(a(j),u)=\{\ell\} \mbox{ for some }u\in G_0(I)\}$.
%
The condition $\{x,y\}\notin\Delta_a$ 
is equivalent to the condition that
$L(a(\{x,y,j\}),u)=\emptyset$ for some $u\in G(I)$. This can also be
refined to the condition that $L(a(j),u)=\{x,y\}$ for some $u\in
G_0(I)$. 
In fact, first of all we have $L(a(\{x,y,j\}),X_i^{\rho_i})=\{i\}\ne
\emptyset$, for all $i=m+1,\ldots, n$, 
since $x,y,j\in [m]$ and $a_i\leq \rho_i-1$. Thus we 
can replace '$u\in G(I)$' in the above condition by 
'$u\in G_0(I)$'.
Also, since $x, y\in{\cal V}_a$,
$L(a(\{x, j\}),u)\ne\emptyset$ 
and 
$L(a(\{y, j\}),u)\ne\emptyset$ 
for all $u\in G_0(I)$. Thus 
$L(a(\{x,y,j\}),u)=\emptyset$  is equivalent to 
$L(a(j),u)=\{x,y\}$ as required.
This is the condition $(ii)(2).$

Now we will show that if $b\in\ZZ^n$ is such that $b\leq \rho-1$, $G_b =
\{j\}$, $L(b(j), u)\ne\emptyset$ for all $u\in G_0(I)$ and
$b(j)\leq a(j)$, then $\Delta_b$ is also connected.  
We prove the contrapositon: if $\Delta_b$ is disjoint
then $\Delta_a$ is disjoint too.
Assume that there exist disjoint non-empty subsets
$P, Q\subset {\cal V}_b$ such that $P\cup Q = {\cal V}_b$ 
and for all $x\in P$ and all $y\in Q$ 
we have $L(b(j), u) = \{x,y\}$ for some $u\in G_0(I)$.
First of all, for $u\in G_0(I)$ we have 
$L(b(\{j, \ell\}),u)\supseteq L(a\{j, \ell\},u)$ so that 
$L_b \subseteq L_a$ and thus ${\cal V}_a \subseteq {\cal V}_b$.
Now for all $x\in P\cap {\cal V}_a$ and $y\in Q\cap{\cal V}_a$,
$L(a(j),u)\subset L(b(j), u)= \{x,y\}$ for some $u\in G_0(I)$.
Also since 
$x, y\in {\cal V}_a$ we must have $L(a(\{j,x\}),u)\ne\emptyset$
and $L(a(\{j,y\}),u)\ne\emptyset$. Then we know that we must 
have $L(a(j), u)=\{x,y\}$. Thus,
by setting $P' = P\cap {\cal V}_a$ and $Q'=Q\cap {\cal V}_a$,
we obtain the non-empty disjoint subsets $P', Q'\subset {\cal V}_a$
showing the disjointness of $\Delta_a$.
Consequently, we can assume $a(j)$ to be maximal as in the
statement $(ii)(c)$. 
Finally, by the proof of Theorem~\ref{finlen-main}, we know that 
the condition $(\ref{4})$ is equivalent to $(iii)$.
\end{proof}

\begin{Remark}
{\em 
Unfortunately we do not know a good combinatorial
characterization for $\tilde{H}_{j}(\Delta_a; K)=0$ for $j\geq 1$,
which is needed to obtain similar results to
Theorem~\ref{finlen-3dim} for $\dim R \geq 4$.
}
\end{Remark}

\subsection{Construction from Buchsbaum Stanley-Reisner ideals}

In this subsection, we compare local cohomologies of 
monomial ideals $I\subset S$ and $\sqrt{I}$.
It is well known that $\Coh{i}{S/\sqrt{I}}_a=0$ for all 
$a\in\ZZ^n$ with $H_a\ne\emptyset$,
which is an immediate consequence from the original version of 
Hochster's formula for Stanley-Reisner ideals.
On the other hand, we may have $\Coh{i}{S/I}_a\ne 0$ for 
such $a\in\ZZ^n$. But for  multi-degrees $a\in\ZZ^n$
with $H_a=\emptyset$, we have an isomorphism.

\begin{Proposition}[Herzog-Takayama-Terai \cite{HTT}]
\label{htt:compare}
Let $I\subset S$ be a monomial ideal. Then we have the following
isomorphisms of $K$-vector spaces
\begin{equation*}
   \Coh{i}{S/I}_a \iso \Coh{i}{S/\sqrt{I}}_a
\end{equation*}
for all $a\in\ZZ^n$ with $H_a=\emptyset$.
\end{Proposition}

\begin{Proposition}[Herzog-Takayama-Terai \cite{HTT}]
\label{htt:genCM}
Let $I\subset S$ be a monomial ideal. Then 
$\length(\Coh{i}{S/I})<\infty$ implies 
$\length(\Coh{i}{S/\sqrt{I}})<\infty$.
In particular, if $I$ is generalized CM, then
$\sqrt{I}$ is also generalized CM (Buchsbaum).
\end{Proposition}
\begin{proof}
We will give here a new proof, which  is
different from that in \cite{HTT}.
Assume that $\length(\Coh{i}{S/I})<\infty$. Then by 
Proposition~\ref{hochster:propFinLen}
$\Coh{i}{S/I}_a=0$ for all $a\in\ZZ^n$ with $G_a\ne\emptyset$.
Thus if $\Coh{i}{S/I}_a\ne 0$ and $H_a=\emptyset$,
we must have $a = (0,\ldots, 0)$.
Now 
by Proposition~\ref{htt:compare} we have 
$\Coh{i}{S/\sqrt{I}}_a\iso \Coh{i}{S/I}_a$ for all $a\in\ZZ^n$
with $H_a=\emptyset$, and this is non-zero
if and only if $a= (0,\ldots, 0)$. Thus we have 
$\length(\Coh{i}{S/\sqrt{I}})<\infty$.
\end{proof}

\begin{Corollary}
Let $I\subset S$ be a generalized CM monimial ideal
and assume that $\sqrt{I}$ is not Cohen-Macaulay (but generalized
 CM by Proposition~\ref{htt:genCM}). Then 
$I$ is not Cohen-Macaulay.
\end{Corollary}
\begin{proof}
Assume that $I$ is Cohen-Macaulay. Then, by 
Proposition~\ref{htt:compare}
and the comment before 
Proposition~\ref{htt:compare},
we have $\Coh{i}{S/\sqrt{I}}_a=0$ for all $i<\dim S/\sqrt{I}$
and for all $a\in\ZZ^n$, namely $\Coh{i}{S/\sqrt{I}}=0$ for all
$i<\dim S/\sqrt{I}$ and $\sqrt{I}$ is Cohen-Macaulay.
\end{proof}

These results suggests a method for constructing (non-CM) generalized
CM monomial ideals from Buchsbaum Stanley-Reisner ideals:
given a Buchsbaum Stanley-Reisner ideal $J$, make 
monomials $X_{j_1}^{e_1}\cdots X_{j_p}^{e_p}$, $(e_i\geq 1, i=1,\ldots,
n)$, for each generator $X_{j_1}\cdots X_{j_p}\in G(J)$. 
In general, one can make more than one monomial generators 
from a single square-free generator.
If we choose suitable exponents $e_i$, the ideal generated by the 
monomials is (non-CM) generalized CM. 
Theorem~\ref{finlen-main}, Corollary~\ref{2dim} and 
Theorem~\ref{finlen-3dim} give the criteria for suitable 
exponents.

\begin{Example}
{\em 
\label{onlyFrobenius}
Let $J = (X_1,\ldots, X_n)(X_{n+1},\ldots, X_{2n})\subset S = \poly{K}{X}{2n}$, 
$(n\geq 2)$.
It is easy to check that $S/I$ is Buchsbaum of dimension $n$ and depth
$1$. Let $I = (X_i^{\alpha_{ij}}X_j^{\beta_{ij}}\;\vert\; 1\leq i \leq n,\; n+1\leq j\leq 2n)$ for 
some $\alpha_{ij}, \beta_{ij}\in\NN$. Then $I$ is generalized CM ideal 
if and only if 
$\alpha_{i,n+1} = \cdots = \alpha_{i,2n}$
for all  $1\leq i \leq n$ 
and 
$\beta_{1,j} = \cdots = \beta_{n, j}$
for all $n+1\leq j\leq 2n$, 
namely $I$ is an image of Frobenius map in the sense of 
Example~\ref{frobenius}.
}
\end{Example}
\begin{Remark}
{\em 
Notice that if $n=1$ then both $J$ and
$I=(X_1^{\alpha_{11}}X_2^{\beta_{11}})$
are Cohen-Macaulay for all $\alpha_{11}, \beta_{11}\in\NN$.
}
\end{Remark}
\begin{proof}[Proof of Example~\ref{onlyFrobenius}]
Assume in the following
that $I$ is generalized CM. Then $I$ must satisfy the 
conditions in Theorem~\ref{finlen-main}, in particular 
the condition for $\length(\Coh{n-1}{S/I}) < \infty$.

Now we set $\rho_i = \max\{\alpha_{ij}\;\vert\; n+1\leq j\leq 2n\}$
and  $\epsilon_i = \min\{\alpha_{ij}\;\vert\; n+1\leq j\leq 2n\}$
for $1\leq i\leq n$ and $\rho_j = \max\{\beta_{ij}\;\vert\; 1\leq i\leq n\}$
and $\epsilon_j = \min\{\beta_{ij}\;\vert\; 1\leq i\leq n\}$ 
for $n+1\leq j\leq 2n$. 
Notice that $\rho_k(\geq 1)$ and $\epsilon_k(\geq 1)$ denote the 
maximal and the minimal exponents of the variable $X_k$, $k=1,\ldots, 2n$.
Let $\Delta$ be the simplicial complex corresponding to 
$J$, which is the disjoint union of two $(n-1)$-simplices
over the vertex set $[n]$ and $[2n]\backslash [n]$. Then
an $(n-2)$-face $\sigma\in\Delta$ is 
either $\sigma = \{1,\ldots, n\}\backslash \{k\}$ for $k=1,\ldots, n$
or $\sigma = \{n+1,\ldots, 2n\}\backslash \{k\}$ for $k=n+1,\ldots, 2n$.

We know that the condition for $a\in\ZZ^{2n}$ such that $0\leq a\leq \rho -1$ 
to be $L(a(\sigma), u)\ne\emptyset$ for all $u\in G(I) = G_0(I)$ is as follows:
\begin{description}
\item [Case (1) $\sigma = \{1,\ldots, n\}\backslash \{k\}$ for $1\leq k\leq n$]
Since 
\begin{eqnarray*}
\lefteqn{
L(a(\sigma), X_i^{\alpha_{ij}}X_j^{\beta_{ij}})
}\\
&=& \left\{
     \begin{array}{ll}
      \{k,j\}   & \mbox{if $i=k$, $a_k<\alpha_{kj}, a_j<\beta_{kj}$} \\
      \{j\}
     & \mbox{if $i=k$, $a_k\geq\alpha_{kj}, a_j<\beta_{kj}$,
                        or if $i\ne k, a_j<\beta_{ij}$} \\
      \{k\}   & \mbox{if $i=k$, $a_k<\alpha_{kj}, a_j\geq \beta_{kj}$} \\
      \emptyset & \mbox{otherwise,}
     \end{array}
    \right.
\end{eqnarray*}
we must have 
\begin{enumerate}
\item $a_j < \min\{\beta_{ij}\;\vert\; 1\leq i\leq n, i\ne k\}$
      for all $n+1\leq j\leq 2n$,
and 
\item for every $n+1\leq j\leq 2n$ we have 
at least one of the followings: (a) $a_k < \alpha_{kj}$, (b) $a_j < \beta_{kj}$.
\end{enumerate}
\item [Case (2) $\sigma = \{n+1,\ldots, 2n\}\backslash \{k\}$ for $n+1\leq k\leq 2n$]
Since 
\begin{eqnarray*}
\lefteqn{
L(a(\sigma), X_i^{\alpha_{ij}}X_j^{\beta_{ij}})
}\\
&=& \left\{
     \begin{array}{ll}
      \{k,i\}   & \mbox{if $j=k$, $a_k<\beta_{ik}, a_i<\alpha_{ik}$} \\
      \{i\}   & \mbox{if $j=k$, $a_k\geq\beta_{ik}, a_i< \alpha_{ik}$}, 
                \mbox{ or if $j\ne k, a_i < \alpha_{ij}$ }\\
      \{k\}     & \mbox{if $j=k$, $a_k<\beta_{ik}, a_i\geq \alpha_{ik}$}\\
      \emptyset & \mbox{otherwise,}
     \end{array}
    \right.
\end{eqnarray*}
we must have 
\begin{enumerate}
\item $a_i < \min\{\alpha_{ij}\;\vert\; n+1\leq j\leq 2n, j\ne k\}$ 
      for all $1\leq i\leq n$,
and 
\item for every $1\leq i\leq n$ we have at least one of the followings:
(a) $a_k < \beta_{ik}$, (b) $a_i < \alpha_{ik}$.
\end{enumerate}
\end{description}
According to Theorem~\ref{finlen-main}, 
for a maximal $a(\sigma)$, where $a$ and $\sigma$ are 
as above, there must exist 
an index $\ell\in [2n]\backslash\sigma$ such that 
$(i)$ $a_\ell =\rho_\ell -1$ and 
$(ii)'$ for every $u=X_i^{\alpha_{ij}}X_j^{\beta_{ij}}
\in G(I)$ with $\nu_\ell(u) = \rho_\ell$
we have $L(a(\sigma\cup\{\ell\}), u)\ne \emptyset$.
Moreover, by the proof of the Theorem~\ref{finlen-main}
we know that 
$(ii)'$ can be replaced by
\begin{equation*}
\mbox{$(ii)$ for every $u\in G(I)$ with $\ell\in\supp(u)$
we have $L(a(\sigma\cup\{\ell\}),u)\ne\emptyset$.
}
\end{equation*}
We now consider the condition for $a$, $\alpha_{ij}$ and $\beta_{ij}$
satisfying  $(i)$ and $(ii)$.
\begin{description}
\item [Case (3) $1\leq \ell\leq n$ with $\ell\notin\sigma$] 
By $(ii)$ we must have 
$\emptyset \ne L(a(\sigma\cup \{\ell\}), 
   X_\ell^{\alpha_{\ell j}}X_{j}^{\beta_{\ell j}}) (\subset \{j\})$
for all $n+1\leq j\leq 2n$.
This holds if and only if $j\notin\sigma$ and 
$a_{j}<\beta_{\ell,j}$.
\item [Case (4) $n+1\leq \ell\leq 2n$ with $\ell\notin\sigma$]
Similarly, we have 
and $a_{i}<\alpha_{i,\ell}$
for any $1\leq i\leq n$ with $i\notin\sigma$.
\end{description}
Now in case (1), an $\ell\in[2n]\backslash\sigma$ satisfying $(i)$ and $(ii)$
must be $\ell = k$, $(1\leq k\leq n)$
or $\ell\in\{n+1,\ldots, 2n\}$. 
Assume that $\ell\in\{n+1,\ldots, 2n\}$. Then the condition 
of case (1) must imply the condition of case (4).
Since 
\begin{eqnarray*}
\max\{\beta_{i\ell}\;\vert\; 1\leq i\leq n, i\ne k\}-1
&\leq &  \max\{\beta_{i\ell}\;\vert\; 1\leq i\leq n\}-1 = \rho_\ell -1 = a_\ell\\
& < &  \min\{\beta_{i\ell}\;\vert\; 1\leq i\leq n, i\ne k\},
\end{eqnarray*}
we have 
$\max\{\beta_{i\ell}\;\vert\; 1\leq i\leq n,i\ne k\}
= \min\{\beta_{i\ell}\;\vert\; 1\leq i\leq n, i\ne k\}$
so that $\beta_{i\ell}$ is constant for all $0\leq i\leq n$ with $i\ne k$.
Now if $\beta_{k\ell}> \beta_{i\ell}$ for some, and equivalently all,
 $i(\ne k)$,  then we have 
$\rho_\ell = \max\{\beta_{i\ell}\;\vert\; 1\leq i\leq n\}= \beta_{k\ell}$
so that $\beta_{k\ell} -1 = \rho_\ell -1 
= a_\ell < \min\{\beta_{i\ell}\;\vert\; 1\leq i\leq n, i\ne k\} = \beta_{i\ell}$.
Then we have $\beta_{k\ell}\leq \beta_{i\ell}$, a contradiction. 
Thus we know $\beta_{k\ell}\leq \beta_{i\ell}$ for all $i\ne k$,
and the index $i$  as in the condition of case (4) 
can be at least $i =1,\ldots, n$ with $i\ne k$ but this 
contradicts the condition $i\notin \sigma = \{1,2,\ldots, n\}\backslash \{k\}$.
Consequently, an $\ell\in[2n]\backslash\sigma$ satisfying $(i)$ and $(ii)$
cannot be from $\{n+1,\ldots, 2n\}$ and 
we must have $\ell = k$ with $1\leq k\leq n$. Now 
the condition of {\bf Case (1)} must imply
the condition of {\bf Case (3)}. 
Comparing {\bf Case (1)} 1 with the condition of {\bf Case (3)}, we know 
that we must have 
\begin{equation}
\label{cond1-3}
  \min\{\beta_{ij}\;\vert\; 1\leq i\leq n, i\ne k\} \leq \beta_{kj}
\end{equation}
for every $j$ with  $n+1\leq j\leq 2n$.
Now we consider similarly with the {\bf Case (2)} and obtain the condition
\begin{equation}
\label{cond2-4}
  \min\{\alpha_{ij}\;\vert\; n+1\leq j\leq 2n, j\ne k\} 
\leq \alpha_{ik}
\end{equation}
for every $i$ with $1\leq i\leq n$ .
Finally, the condition $(\ref{cond1-3})$ for $k=1,\ldots, n$
together with the condition $(\ref{cond2-4})$ for $k=n+1,\ldots, 2n$ 
entails $\beta_{ij}$ are constant for all $1\leq i\leq n$ 
and $\alpha_{ij}$ are constant for all $n+1\leq j\leq 2n$, i.e.,
$I$ is obtained from $J$ by Frobenius transformation
in the  sense of Example~\ref{frobenius}.
\end{proof}

\begin{Example}
{\em 
If we allow to make more than two monomial generators from a 
single square-free generator, we can construct 
more generalized CM monomial ideals from the same 
Stanley-Reisner ideals as in Example~\ref{onlyFrobenius}, For example,
from $J_1 = (X_1,X_2)(X_3,X_4)\subset K[X_1,X_2,X_3,X_4]$
we make
\begin{equation*}
I_1 = (X_1X_3, X_1^2X_4, X_1X_4^2, X_2^2X_3, X_2X_3^2, X_2X_4).
\end{equation*}
Also from 
$J_2 = (X_1,X_2,X_3)(X_4,X_5,X_6)
\subset K[X_1, \ldots, X_6]$
we make
\begin{equation*}
I_2 = (X_1^3X_4,X_1X_4^5, X_1X_5, X_1X_6, X_2X_4, X_2X_5, X_2X_6,
 X_3X_4, X_3X_5, X_3X_6).
\end{equation*}
$I_1$ and $I_2$ are both generalized CM, but for example
\begin{equation*}
I_3 = (X_1^3X_4,X_1^2X_4^2, X_1X_4^3, X_1X_5, X_1X_6, X_2X_4, X_2X_5, X_2X_6, X_3X_4, X_3X_5, X_3X_6)
\end{equation*}
is not generalized CM.
}
\end{Example}


\end{document}